\documentclass[10pt]{amsart}

\usepackage[T1]{fontenc}
\usepackage{lmodern}
\usepackage{amsmath,amssymb,mathtools}
\usepackage[expansion=false]{microtype}
\usepackage[hidelinks]{hyperref}

\newcommand{\Q}{\mathbf{Q}}
\newcommand{\Z}{\mathbf{Z}}
\newcommand{\C}{\mathbf{C}}
\newcommand{\R}{\mathbf{R}}
\newcommand{\cO}{\mathcal{O}}
\newcommand{\cA}{\mathcal{A}}
\newcommand{\cD}{\mathcal{D}}
\newcommand{\cX}{\mathcal{X}}
\newcommand{\one}{\mathbf{1}}

\theoremstyle{plain}
\newtheorem{theorem}{Theorem}[section]
\newtheorem{proposition}[theorem]{Proposition}
\newtheorem{lemma}[theorem]{Lemma}
\newtheorem{corollary}[theorem]{Corollary}
\theoremstyle{remark}
\newtheorem{remark}[theorem]{Remark}
\newtheorem{example}[theorem]{Example}

\title[Integral trace forms of cyclic number fields]{Integral trace forms of cyclic number fields without tameness}

\author{Zhi-Lin Zhang}
\address{Independent Researcher, Taipei, Taiwan}
\email{hsa00000@gmail.com}

\subjclass[2020]{Primary 11R04; Secondary 11R20, 11R33, 11E12, 11L05}
\keywords{cyclic number field, integral trace form, discriminant, Dirichlet character, associated order, Gaussian period, wild ramification}

\begin{document}

\begin{abstract}
Let \(K/\Q\) be a cyclic number field and let
\[
T_K=\bigl(\cO_K,(x,y)\mapsto \operatorname{Tr}_{K/\Q}(xy)\bigr)
\]
be its integral trace lattice.  We prove that two cyclic number fields have
isometric integral trace lattices if and only if they have the same degree,
signature, and signed discriminant.  More generally, an isomorphism between
the Dirichlet character groups of two finite abelian extensions that preserves
the conductor and parity of every character yields an isometry of their
integral trace lattices, equivariant for the dual identification of Galois
groups.  For cyclic fields, the degree, signature, and
discriminant determine the character group together with these two functions.
This proves Bola\~nos's conjecture for wildly ramified cyclic fields and gives
the corresponding extension of the tame classification: in general, the
signature must be retained as an additional invariant.
\end{abstract}

\maketitle

\section{Introduction}

For a number field \(K\), the integral trace form is the symmetric integral bilinear form
\[
B_K\colon \cO_K\times \cO_K\longrightarrow \Z,
\qquad
B_K(x,y)=\operatorname{Tr}_{K/\Q}(xy).
\]
We write \(T_K=(\cO_K,B_K)\).  Its rank is \([K:\Q]\), its determinant is the signed field discriminant \(\operatorname{Disc}(K)\), and its real signature determines the signature of \(K\).  Throughout,
\(\operatorname{sig}(K)=(r_1(K),r_2(K))\).  Thus the degree, the discriminant,
and the signature are necessary invariants for the isometry class of \(T_K\).

Conner and Perlis proved that equal discriminants yield isometric integral
trace forms for two cyclic fields of the same odd prime degree \(p\)
\cite[Chapter~IV]{ConnerPerlis}.  Bola\~nos and Mantilla-Soler established the
corresponding classification for two tame cyclic fields of the same degree \(n\)
\cite[Theorems~4.2 and~4.5]{BolanosMantilla}; their explicit odd- and even-degree
Gram matrices appear in
\cite[Corollary~4.4 and Section~4.2]{BolanosMantilla}.  Their proof uses the
Hilbert--Speiser theorem and does not cover general wild ramification.  Bola\~nos
conjectured the analogous statement for wildly ramified cyclic fields of the
same degree, signature, and discriminant
\cite[Section~3.3, p.~72]{BolanosThesis}.  The main result establishes this
conjectured wild-ramification case with fixed signature.

\begin{samepage}
\begin{theorem}[Main theorem]\label{thm:main}
Let \(K\) and \(L\) be cyclic number fields.  Then
\[
T_K\cong T_L
\]
if and only if
\[
[K:\Q]=[L:\Q],
\qquad
\operatorname{sig}(K)=\operatorname{sig}(L),
\qquad
\operatorname{Disc}(K)=\operatorname{Disc}(L).
\]
\end{theorem}
\end{samepage}

The signature hypothesis cannot be dropped in general.  Already in degree
\(4\), the cyclic fields
\[
\Q\bigl(\sqrt{2+\sqrt2}\bigr)
\quad\text{and}\quad
\Q\bigl(\sqrt{-2-\sqrt2}\bigr)
\]
have the same signed discriminant \(2^{11}\), but one is totally real and the
other is totally imaginary.  Their trace forms therefore have different real
signatures.  A full calculation of this example is included below.
Since \(\operatorname{Disc}(K)=(-1)^{r_2(K)}|\operatorname{Disc}(K)|\), the
signed discriminant in the theorem can equivalently be replaced by its
absolute value when the signature is retained.

The proof separates arithmetic reconstruction from integral-lattice comparison.
At each odd prime, the degree together with the discriminant exponent uniquely
recovers the embedded local cyclic character subgroup.  At \(2\), the degree and
discriminant recover the subgroup except for one even--odd ambiguity, which the
global signature resolves.  This produces a character-group isomorphism
preserving the conductor and parity of every character.

The comparison step applies more generally to arbitrary finite abelian
extensions.  It uses Leopoldt's description of the ring of integers as a free
rank-one module over its associated order \cite{Leopoldt,Lettl}.  The associated
order is controlled by character conductors, while the trace pairing is expressed using normalized character resolvents and the
product of the Gauss sums of a character and its inverse.  The resulting criterion shows that the abstract character
group together with the
conductor and parity of every character already determines the integral trace
lattice.  Thus the cyclic classification follows from a more general comparison theorem
for finite abelian extensions.

\section{Character groups and local conductors}\label{sec:prelim}

This section fixes the character-theoretic invariants used to recover local
ramification and to compare trace lattices.

All number fields under consideration are embedded in \(\C\).  A Dirichlet
character will be viewed as a finite-order continuous character of
\(\widehat{\Z}^{\,\times}\).  Its conductor is the least positive integer \(m\)
through which it factors via \((\Z/m\Z)^\times\).  This convention makes
multiplication of characters independent of a choice of modulus; when a
primitive Dirichlet character is needed, we use the character at its
conductor.

Let \(E/\Q\) be a finite abelian extension, let
\(G_E=\operatorname{Gal}(E/\Q)\), and let
\[
X_E=\widehat{G_E}=\operatorname{Hom}(G_E,\C^\times).
\]
By Kronecker--Weber, pullback along the cyclotomic Galois action identifies
\(X_E\) with a group of Dirichlet characters.  We denote the conductor of
\(\chi\in X_E\) by \(f_\chi\), with \(f_{\one}=1\).  The conductor of the field is
\[
f_E=\operatorname{lcm}_{\chi\in X_E}f_\chi.
\]
It is the least \(f\) such that \(E\subseteq\Q(\zeta_f)\), and it is never
congruent to \(2\) modulo \(4\).
We will use the character group \(X_E\) together with the two functions
\[
\chi\longmapsto f_\chi,
\qquad
\chi\longmapsto\chi(-1).
\]
Accordingly, when comparing two such character groups, we will explicitly
require the group isomorphism to preserve both conductor and parity.

The conductor--discriminant formula states that
\begin{equation}\label{eq:conductor-discriminant}
|\operatorname{Disc}(E)|=\prod_{\chi\in X_E}f_\chi.
\end{equation}
We use the standard character-theoretic form of this formula; see, for example,
\cite[Theorem~3.11]{Washington}.

Every Dirichlet character has a unique factorization
\[
\chi=\prod_p\chi_p,
\]
where \(\chi_p\) has \(p\)-power conductor.  Let \(\cX_p\) denote the group of all finite-order Dirichlet characters of \(p\)-power conductor, and let
\[
\pi_p\colon X_E\longrightarrow \cX_p,
\qquad
\chi\longmapsto\chi_p
\]
be the local projection.  Put
\[
Y_{E,p}=\pi_p(X_E).
\]
The kernel of \(\pi_p\) consists of the characters trivial on inertia at \(p\).
Character duality therefore identifies \(Y_{E,p}\) with the dual of the inertia
group, so \(|Y_{E,p}|\) is the ramification index of \(p\) in \(E/\Q\).

For \(\eta\in\cX_p\), define
\[
c_p(\eta)=v_p(f_\eta).
\]
The next formula is immediate from \eqref{eq:conductor-discriminant}.

\begin{lemma}[Local discriminant formula]\label{lem:localdisc}
Let \(E/\Q\) be abelian of degree \(n\).  Then, for every prime \(p\),
\begin{equation}\label{eq:localdisc}
v_p(|\operatorname{Disc}(E)|)
=
\frac{n}{|Y_{E,p}|}
\sum_{\eta\in Y_{E,p}}c_p(\eta).
\end{equation}
\end{lemma}

\begin{proof}
The surjection \(\pi_p\colon X_E\to Y_{E,p}\) has fibres of size \(n/|Y_{E,p}|\).  Since
\[
v_p(f_\chi)=c_p(\chi_p),
\]
summing over \(X_E\) gives \eqref{eq:localdisc}.  
\end{proof}

The trace lattice itself recovers the field signature.

\begin{proposition}\label{prop:necessary}
Let \(E\) have signature \((r_1,r_2)\).  The real bilinear form \(T_E\otimes_\Z\R\) has signature
\[
(r_1+r_2,r_2).
\]
Consequently, an isometry \(T_E\cong T_F\) implies equality of degrees, field signatures, and signed discriminants.
\end{proposition}

\begin{proof}
Under the Minkowski decomposition
\[
E\otimes_\Q\R\cong \R^{r_1}\times\C^{r_2},
\]
a real embedding contributes the positive form \(x\mapsto x^2\).  On a complex factor, writing \(z=a+bi\), one has
\[
\operatorname{Tr}_{\C/\R}(z^2)=2(a^2-b^2),
\]
which has signature \((1,1)\).  This gives the asserted signature.  The rank is the degree, and the determinant of a Gram matrix of \(T_E\) is \(\operatorname{Disc}(E)\).  
\end{proof}

\section{Local character groups at odd primes}\label{sec:odd}

We first show that a local character subgroup at an odd prime can be
reconstructed from its contribution to the discriminant.

Fix an odd prime \(p\).  The group of finite-order characters of \(p\)-power conductor has the abstract structure
\begin{equation}\label{eq:oddlocalgroup}
\cX_p\cong C_{p-1}\times C_{p^\infty},
\end{equation}
where \(C_{p^\infty}\) denotes the Pr\"ufer \(p\)-group.  The Sylow subgroups
of every finite subgroup are cyclic and have pairwise coprime orders; hence
every finite subgroup is cyclic.  In particular, for every \(d\mid p-1\) and
\(s\ge0\), there is a unique subgroup of order \(dp^s\).

For \(s\ge0\), put
\begin{equation}\label{eq:Sps}
S_{p,s}=\sum_{j=1}^{s}\varphi(p^j)(j+1),
\qquad S_{p,0}=0.
\end{equation}

\begin{lemma}\label{lem:oddconductorsum}
Let \(Y\subset\cX_p\) be the cyclic subgroup of order \(dp^s\), where \(d\mid p-1\).  Then
\begin{equation}\label{eq:oddconductorsum}
\sum_{\eta\in Y}c_p(\eta)
=(d-1)+dS_{p,s}.
\end{equation}
\end{lemma}

\begin{proof}
The elements with trivial \(p\)-primary component form the unique subgroup of order \(d\).  Its identity has conductor exponent \(0\), while its other \(d-1\) elements are nontrivial tame characters and have conductor exponent \(1\).

For \(1\le j\le s\), there are \(d\varphi(p^j)\) elements whose \(p\)-primary component has exact order \(p^j\).  Each such character has conductor exponent \(j+1\).  Summing these contributions proves \eqref{eq:oddconductorsum}.  
\end{proof}

Combining Lemmas~\ref{lem:localdisc} and \ref{lem:oddconductorsum}, if \(E/\Q\) is abelian of degree \(n\) and \(|Y_{E,p}|=dp^s\), then
\begin{equation}\label{eq:oddap}
\frac{v_p(|\operatorname{Disc}(E)|)}{n}
=
\frac{1}{p^s}
\left(S_{p,s}+1-\frac1d\right).
\end{equation}
We next show that the possible values for distinct \(s\) are separated; for fixed \(s\), the right-hand side is strictly increasing in \(d\). Consequently, it uniquely determines both \(s\) and \(d\).

\begin{lemma}[Separation for distinct values of \(s\)]\label{lem:sseparation}
For every odd prime \(p\) and every \(s\ge0\),
\begin{equation}\label{eq:gapidentity}
\frac{S_{p,s+1}}{p^{s+1}}
-
\frac{S_{p,s}+1}{p^s}
=
\frac{p^{s+1}-2}{p^{s+1}}>0.
\end{equation}
\end{lemma}

\begin{proof}
Writing \(\varphi(p^j)=p^j-p^{j-1}\) and telescoping gives
\[
S_{p,s}=(s+1)p^s-1-\frac{p^s-1}{p-1},
\]
including \(s=0\).  Substitution for \(s\) and \(s+1\) gives
\eqref{eq:gapidentity}.  Its right side is positive because \(p\ge3\).  
\end{proof}

\begin{proposition}[Determination of \(Y_{E,p}\) for odd \(p\)]\label{prop:odddetermination}
Let \(E/\Q\) be abelian of degree \(n\), and let \(p\) be odd.  Then \(n\) and
\(v_p(|\operatorname{Disc}(E)|)\) uniquely determine the embedded subgroup
\(Y_{E,p}\subset\cX_p\).
\end{proposition}

\begin{proof}
For fixed \(s\), the right side of \eqref{eq:oddap} is strictly increasing as a function of \(d\).  Moreover, as \(d\) varies, it lies in the half-open interval
\[
\left[
\frac{S_{p,s}}{p^s},
\frac{S_{p,s}+1}{p^s}
\right).
\]
Lemma~\ref{lem:sseparation} shows that these intervals are disjoint for distinct \(s\).  Hence the discriminant exponent first determines \(s\) and then determines \(d\).  The uniqueness in \eqref{eq:oddlocalgroup} then determines \(Y_{E,p}\) itself.  
\end{proof}

\begin{remark}\label{rem:oddparity}
If \(Y\subset\cX_p\) is cyclic, all generators of \(Y\) have the same value at \(-1\).  Indeed, if \(|Y|\) is odd, every character in \(Y\) is even.  If \(|Y|\) is even, two generators differ by an exponent prime to \(|Y|\), hence by an odd exponent, which preserves the value \(\pm1\) at \(-1\).  Thus Proposition~\ref{prop:odddetermination} also determines the parity of a generator of \(Y_{E,p}\).
\end{remark}

\section{The dyadic local character group}\label{sec:two}

We next isolate the one ambiguity that the dyadic discriminant contribution
does not detect and show that it is exactly a parity ambiguity.

For \(a\ge3\),
\[
(\Z/2^a\Z)^\times\cong C_2\times C_{2^{a-2}}.
\]
For every \(m\ge1\), fix \(\zeta_m=\exp(2\pi i/m)\).
Let \(\varepsilon\) be the primitive character of conductor \(4\) with
\[
\varepsilon(-1)=-1,
\qquad
\varepsilon(5)=1.
\]
For \(s\ge1\), choose the even character \(\psi_s\) characterized by
\[
\psi_s(-1)=1,
\qquad
\psi_s(5)=\zeta_{2^s}.
\]
Here \(\psi_s\) is regarded as a character modulo \(2^{s+2}\) and then inflated
to \(\widehat{\Z}^{\,\times}\).  It has exact order \(2^s\) and conductor
\(2^{s+2}\), and \(\psi_s^2=\psi_{s-1}\) for \(s\ge2\).
The dyadic character group has structure
\[
\cX_2\cong C_2\times C_{2^\infty},
\]
where the first factor is generated by \(\varepsilon\) and the second consists
of the even characters.

\begin{lemma}[Cyclic dyadic subgroups]\label{lem:dyadicsubgroups}
The nontrivial finite cyclic subgroups of \(\cX_2\) are as follows:
\begin{enumerate}
\item the three order-\(2\) subgroups
\[
\langle\varepsilon\rangle,
\qquad
\langle\psi_1\rangle,
\qquad
\langle\varepsilon\psi_1\rangle;
\]
\item for each \(s\ge2\), the two order-\(2^s\) subgroups
\[
\langle\psi_s\rangle
\qquad\text{and}\qquad
\langle\varepsilon\psi_s\rangle.
\]
The two subgroups in item~(2) are distinguished by the parity of a generator.
\end{enumerate}
\end{lemma}

\begin{proof}
The subgroup of characters of order dividing \(2^s\) is
\(\langle\varepsilon,\psi_s\rangle\cong C_2\times C_{2^s}\).
For \(s\ge2\), an element of exact order \(2^s\) has the form
\(\varepsilon^a\psi_s^u\) with \(a\in\{0,1\}\) and \(u\) odd.
Raising it to an exponent inverse to \(u\) modulo \(2^s\) changes \(u\) to \(1\)
and leaves \(a\) unchanged.  The resulting two subgroups are distinct: the
first contains only even characters, whereas a generator of the second is
odd.  When \(s=1\), the subgroup of characters of order dividing \(2\) is
\(C_2\times C_2\) and has the three stated order-\(2\) subgroups.  
\end{proof}

\begin{lemma}[Dyadic conductor sums]\label{lem:dyadicsums}
Let \(Y\subset\cX_2\) be cyclic.
\begin{enumerate}
\item If \(Y=\langle\varepsilon\rangle\), then
\[
\sum_{\eta\in Y}c_2(\eta)=2.
\]
\item If \(Y=\langle\psi_1\rangle\) or \(Y=\langle\varepsilon\psi_1\rangle\), then
\[
\sum_{\eta\in Y}c_2(\eta)=3.
\]
\item If \(|Y|=2^s\) with \(s\ge2\), then
\begin{equation}\label{eq:dyadicsum}
\sum_{\eta\in Y}c_2(\eta)
=3+\sum_{j=2}^{s}\varphi(2^j)(j+2)
=(s+1)2^s-1.
\end{equation}
\end{enumerate}
\end{lemma}

\begin{proof}
The first two assertions follow from the conductors \(4\) and \(8\).  For \(s\ge2\), the unique character of order \(2\) in either subgroup is \(\psi_1\), of conductor exponent \(3\).  For each \(2\le j\le s\), there are \(\varphi(2^j)\) characters of exact order \(2^j\), each of conductor exponent \(j+2\).  This proves the first equality in \eqref{eq:dyadicsum}; the second follows by induction on \(s\).  
\end{proof}

Now let \(E/\Q\) be cyclic of degree \(n\).  Then \(Y_{E,2}\) is cyclic, so
Lemma~\ref{lem:localdisc} gives the following possible dyadic discriminant
exponents:
\begin{equation}\label{eq:dyadicexponents}
v_2(|\operatorname{Disc}(E)|)=
\begin{cases}
0,&Y_{E,2}=1,\\
n,&Y_{E,2}=\langle\varepsilon\rangle,\\
\dfrac{3n}{2},&Y_{E,2}=\langle\psi_1\rangle
\text{ or }\langle\varepsilon\psi_1\rangle,\\
n(s+1)-\dfrac{n}{2^s},&|Y_{E,2}|=2^s,\ s\ge2.
\end{cases}
\end{equation}
For the exponent \(3n/2\), and for each exponent on the last line, the two
possible subgroups have the same order and conductor distribution but opposite
generator parity.

\begin{proposition}[Determination of \(Y_{E,2}\) up to generator parity]\label{prop:dyadicdetermination}
Let \(E/\Q\) be cyclic of degree \(n\).  The values \(n\) and
\(v_2(|\operatorname{Disc}(E)|)\) determine \(Y_{E,2}\) except in the following
cases:
\begin{enumerate}
\item if the exponent is \(3n/2\), then
\(Y_{E,2}=\langle\psi_1\rangle\) or
\(Y_{E,2}=\langle\varepsilon\psi_1\rangle\);
\item if the exponent is
\(A_s=n(s+1)-n/2^s\) for some \(s\ge2\), then
\(Y_{E,2}=\langle\psi_s\rangle\) or
\(Y_{E,2}=\langle\varepsilon\psi_s\rangle\).
\end{enumerate}
In each case the two possibilities are distinguished by the parity of a
generator.
\end{proposition}

\begin{proof}
The first three values in \eqref{eq:dyadicexponents} are \(0\), \(n\), and \(3n/2\).  For \(s\ge2\), put
\[
A_s=n(s+1)-\frac{n}{2^s}.
\]
Then \(A_2=11n/4>3n/2\) and
\[
A_{s+1}-A_s=n+\frac{n}{2^{s+1}}>0.
\]
Thus the exponent determines either the unique subgroup on its line or one of
the two subgroups listed in the proposition.  Lemma~\ref{lem:dyadicsubgroups}
distinguishes each ambiguous pair by generator parity.  
\end{proof}

To resolve that ambiguity globally, we use the relation between parity and
signature.  Complex conjugation corresponds, under character duality, to
evaluation at \(-1\).

\begin{lemma}[Parity and signature]\label{lem:paritysignature}
Let \(E/\Q\) be abelian.
\begin{enumerate}
\item \(E\) is totally real if and only if \(\chi(-1)=1\) for every \(\chi\in X_E\).
\item If \(X_E\) is cyclic with generator \(\chi\), then
\[
\chi(-1)=
\begin{cases}
1,&E\text{ is totally real},\\
-1,&E\text{ is totally imaginary}.
\end{cases}
\]
\end{enumerate}
\end{lemma}

\begin{proof}
The fixed field of complex conjugation is the maximal totally real subfield.
Thus complex conjugation acts trivially on \(E\) exactly when every character of
\(G_E\) is trivial on it, which is equivalent to \(\chi(-1)=1\) for every
\(\chi\in X_E\).  Since \(E/\Q\) is Galois, either all embeddings are real or none is real.
If \(X_E\) is cyclic, it is enough to test a generator; a nontrivial value on
complex conjugation must be \(-1\).  
\end{proof}

For cyclic global fields, the signature resolves this ambiguity.

\begin{proposition}[Determination of all local character subgroups]\label{prop:localdetermination}
Let \(K/\Q\) be cyclic of degree \(n\).  Then the triple
\[
\bigl(n,|\operatorname{Disc}(K)|,\operatorname{sig}(K)\bigr)
\]
uniquely determines \(Y_{K,p}\subset\cX_p\) for every prime \(p\).
\end{proposition}

\begin{proof}
For odd \(p\), this is Proposition~\ref{prop:odddetermination}.  Let \(X_K=\langle\chi\rangle\) and write \(\chi=\prod_p\chi_p\).  Each \(\chi_p\) is a generator of \(Y_{K,p}\).  By Remark~\ref{rem:oddparity}, the values \(\chi_p(-1)\) for odd \(p\) are already determined by the discriminant.  Lemma~\ref{lem:paritysignature} determines \(\chi(-1)\) from the signature, and therefore
\begin{equation}\label{eq:dyadicparity}
\chi_2(-1)=\chi(-1)\prod_{p\text{ odd}}\chi_p(-1).
\end{equation}
Proposition~\ref{prop:dyadicdetermination} and \eqref{eq:dyadicparity} determine \(Y_{K,2}\) completely.  
\end{proof}

\begin{theorem}[Conductor- and parity-preserving character-group isomorphism]\label{thm:characterdata}
Let \(K\) and \(L\) be cyclic number fields with the same degree, absolute discriminant, and signature.  Then there is a group isomorphism
\[
\iota\colon X_K\xrightarrow{\sim}X_L
\]
such that
\[
f_{\iota(\chi)}=f_\chi
\qquad\text{and}\qquad
\iota(\chi)(-1)=\chi(-1)
\]
for every \(\chi\in X_K\).
\end{theorem}

\begin{proof}
Choose generators \(X_K=\langle\chi_K\rangle\) and \(X_L=\langle\chi_L\rangle\).  Proposition~\ref{prop:localdetermination} gives
\[
Y_{K,p}=Y_{L,p}
\qquad\text{for every }p.
\]
Hence, for every \(p\), there is an integer \(u_p\) relatively prime to
\(|Y_{K,p}|\) such that
\[
(\chi_L)_p=(\chi_K)_p^{u_p};
\]
when \(Y_{K,p}=1\), take \(u_p=1\).  Define
\[
\iota\colon X_K\longrightarrow X_L,
\qquad \iota(\chi_K^j)=\chi_L^j.
\]
The assignment is well-defined because both generators have the same order
\(n\), and it is a group isomorphism.  For any finite-order character \(\eta\),
raising \(\eta\) to a power relatively prime to its order leaves its kernel
unchanged and hence leaves its conductor unchanged.  Applying this to
\(\eta=(\chi_K)_p^j\) gives
\[
v_p(f_{\chi_K^j})=v_p(f_{\chi_L^j})
\]
for every \(p\) and every integer \(j\).  Thus the global conductors agree.

By the equality of signatures and Lemma~\ref{lem:paritysignature},
\[
\chi_K(-1)=\chi_L(-1).
\]
Raising to the \(j\)th power proves preservation of parity for every character.  
\end{proof}

\section{Leopoldt's associated-order decomposition}\label{sec:leopoldt}

We now prepare the integral comparison by expressing the ring of integers using subsets of characters determined by their
conductors.

Let \(E/\Q\) be an arbitrary finite abelian extension with Galois group \(G\),
character group \(X=\widehat G\), and conductor \(f=f_E\).  For a positive integer
\(m\), define
\[
p(m)=\prod_{\substack{p\mid m\\p\ne2}}p,
\qquad
q(m)=\prod_{v_p(m)\ge2}p^{v_p(m)}.
\]
Put
\[
\cD(f)=
\{m\ge1:p(f)\mid m,\ m\mid f,\ m\not\equiv2\pmod4\}.
\]
For \(m\in\cD(f)\), put
\[
\Phi_m=\{\chi\in X:q(f_\chi)=q(m)\}.
\]
These sets give a disjoint union
\(X=\coprod_{m\in\cD(f)}\Phi_m\); some of them may be empty.
For each \(\chi\in X\), the unique \(m\in\cD(f)\) such that
\(\chi\in\Phi_m\) is
\begin{equation}\label{eq:mchi}
m(\chi)=f_\chi
\prod_{\substack{p\mid f,\ p\nmid f_\chi\\p\ne2}}p.
\end{equation}
Indeed, every odd prime dividing \(f\) must occur in \(m(\chi)\), an exponent
at least \(2\) must agree with that in \(f_\chi\), and an exponent \(1\) at \(2\)
is excluded.  These requirements give exactly the displayed integer.

For \(\chi\in X\), let
\[
e_\chi=\frac1{|G|}\sum_{g\in G}\chi(g)^{-1}g\in\C[G],
\]
and for \(m\in\cD(f)\) put
\[
e_m=\sum_{\chi\in\Phi_m}e_\chi.
\]
The conductor is unchanged by a Galois conjugation of character values:
such a conjugation replaces a character by a power relatively prime to its
order and preserves its kernel.  Thus \(\Phi_m\) is stable under every
Galois conjugation, so \(e_m\in\Q[G]\).  The \(e_m\) are pairwise orthogonal
idempotents with sum \(1\); if \(\Phi_m\) is empty, then \(e_m=0\).
The group algebra \(\Q[G]\) acts on \(E\) by
\[
\left(\sum_{g\in G}a_g g\right)x=\sum_{g\in G}a_g g(x).
\]

Let
\[
E_m=E\cap\Q(\zeta_m),
\qquad
d_{E,m}=[E_m:\Q],
\]
and define the Gaussian period
\[
\eta_{E,m}
=
\operatorname{Tr}_{\Q(\zeta_m)/E_m}(\zeta_m).
\]
Finally set
\[
\Theta_E=\sum_{m\in\cD(f)}\eta_{E,m}.
\]
The associated order is
\[
\cA_E=
\{a\in\Q[G]:a\cO_E\subseteq\cO_E\}.
\]

The intersection degrees are determined by the character conductors.

\begin{lemma}\label{lem:dmcharacters}
Let \(E/\Q\) be a finite abelian extension with character group \(X_E\) and
conductor \(f_E\).  For \(m\mid f_E\),
\begin{equation}\label{eq:dmcharacters}
d_{E,m}
=
[E\cap\Q(\zeta_m):\Q]
=
\#\{\chi\in X_E:f_\chi\mid m\}.
\end{equation}
\end{lemma}

\begin{proof}
The character group of \(E\cap\Q(\zeta_m)\) consists exactly of those characters
of \(E\) that factor through \((\Z/m\Z)^\times\), equivalently those whose
conductor divides \(m\).  Its cardinality is the degree of the intersection.

\end{proof}

We use the following theorem of Leopoldt, in the form given by Lettl
\cite{Leopoldt,Lettl}; a modern restatement appears in
\cite[Theorem~3]{Aoki}.

\begin{theorem}[Leopoldt--Lettl]\label{thm:leopoldt}
With the notation above,
\begin{equation}\label{eq:leopoldtdecomp}
\cO_E
=
\bigoplus_{m\in\cD(f)}\Z[G]\eta_{E,m}
=
\cA_E\Theta_E,
\end{equation}
and \(\cO_E\) is free of rank one over \(\cA_E\).  Moreover,
\begin{equation}\label{eq:associatedorder}
\cA_E=\Z[G][e_m:m\in\cD(f)].
\end{equation}
\end{theorem}

\section{Gaussian periods and character resolvents}\label{sec:gauss}

This section computes the character resolvents of the element \(\Theta_E\)
and uses Fourier inversion to recover the trace pairing.  The fixed embedding
\(E\subseteq\C\) is used whenever an element of \(E\) occurs in a complex sum.

For a Dirichlet character \(\chi\) whose conductor divides \(m\), define
\[
\tau_m(\chi)=
\sum_{a\in(\Z/m\Z)^\times}\chi(a)\zeta_m^a,
\]
where \(\chi\) in this sum is the character induced modulo \(m\).
Set \(\tau_1(\one)=1\).  Let \(\mu\) denote the M\"obius function.
For a primitive character \(\chi^*\) of conductor \(d\), we use the standard
identity
\begin{equation}\label{eq:primitivegauss}
\tau_d(\chi^*)\tau_d((\chi^*)^{-1})=\chi^*(-1)d
\end{equation}
from \cite[Chapter~4]{Washington}.  In particular, a primitive Gauss sum is
nonzero.  Our convention makes \eqref{eq:primitivegauss} valid at \(d=1\).

\begin{lemma}[Gauss sums at an induced modulus]\label{lem:imprimitivegauss}
Let \(\chi\) have conductor \(d\mid m\), let \(\chi^*\) be its primitive
character modulo \(d\), and put \(r=m/d\).  Extend \(\chi^*\) to integers by zero
on integers not relatively prime to \(d\).  Then
\begin{equation}\label{eq:imprimitivegauss}
\tau_m(\chi)=\mu(r)\chi^*(r)\tau_d(\chi^*).
\end{equation}
In particular, \(\tau_m(\chi)\ne0\) if and only if \(r\) is squarefree and
\(\gcd(r,d)=1\).
\end{lemma}

\begin{proof}
Let \(P\) be the product of the primes dividing \(m\) but not \(d\).
The value \(\chi^*(a)\) already vanishes when \(a\) is not relatively prime to
\(d\).  Inclusion--exclusion for the remaining primes gives
\[
\tau_m(\chi)=
\sum_{e\mid P}\mu(e)\chi^*(e)
\sum_{b\bmod m/e}\chi^*(b)\zeta_{m/e}^{\,b}.
\]
Since \(\gcd(e,d)=1\), we have \(d\mid m/e\).
Write \(m/e=dq\) and \(b=b_0+dt\) with \(b_0\bmod d\) and \(0\le t<q\).
The sum over \(t\) is \(\sum_{t=0}^{q-1}\zeta_q^t\), which vanishes unless
\(q=1\).  Hence the only possible nonzero term is \(e=m/d=r\), and this term
occurs precisely when \(r\mid P\).  This proves \eqref{eq:imprimitivegauss}.
A primitive Gauss sum is nonzero, so the last assertion follows.  For
\(d=1\), the same argument uses the constant character \(\one\) and gives
\(\tau_m(\one)=\mu(m)\).  
\end{proof}

\begin{lemma}[Criterion for membership in \(\Phi_m\)]\label{lem:phimcriterion}
If \(\chi\in\Phi_m\), then \(f_\chi\mid m\) and
\[
r=\frac{m}{f_\chi}
\]
is squarefree and relatively prime to \(f_\chi\).
Conversely, for \(m\in\cD(f)\) with \(f_\chi\mid m\), these two conditions
imply \(\chi\in\Phi_m\).
\end{lemma}

\begin{proof}
The forward implication follows directly from \eqref{eq:mchi}.
For the converse, multiplication by a squarefree integer relatively prime
to \(f_\chi\) does not change its factors of prime-power exponent at least
\(2\).  Thus \(q(m)=q(f_\chi)\), which is the definition of \(\chi\in\Phi_m\).

\end{proof}

\begin{lemma}[Gauss-sum product for inverse characters]\label{lem:gaussproduct}
If \(\chi\in\Phi_m\), then
\begin{equation}\label{eq:gaussproduct}
\tau_m(\chi)\tau_m(\chi^{-1})=\chi(-1)f_\chi.
\end{equation}
\end{lemma}

\begin{proof}
Put \(d=f_\chi\) and \(r=m/d\).  Lemma~\ref{lem:phimcriterion} gives
\(\mu(r)^2=1\) and
\(\chi^*(r)(\chi^*)^{-1}(r)=1\).
Applying \eqref{eq:imprimitivegauss} to \(\chi\) and \(\chi^{-1}\), the extra
factors cancel by these two equalities.  Formula~\eqref{eq:primitivegauss}
then gives the result.  
\end{proof}

For \(n=|G|\) and \(\chi\in X\), define the normalized character resolvent
\begin{equation}\label{eq:resolvent}
\rho_{E,\chi}\colon E\longrightarrow\C,
\qquad
\rho_{E,\chi}(x)=\frac1n\sum_{g\in G}\chi(g)^{-1}g(x).
\end{equation}
This is a \(\Q\)-linear map with values in \(\C\) under the fixed embedding.
For \(a=\sum_g a_g g\in\C[G]\), define its Fourier transform
\(\widehat a\colon X\to\C\) by
\[
\widehat a(\chi)=\sum_{g\in G}a_g\chi(g).
\]
For \(h\in G\), changing the summation variable in \eqref{eq:resolvent} gives
\[
\rho_{E,\chi}(h(x))
=\frac1n\sum_{k\in G}\chi(kh^{-1})^{-1}k(x)
=\chi(h)\rho_{E,\chi}(x).
\]
Consequently, for \(a\in\Q[G]\) and \(x\in E\),
\begin{equation}\label{eq:resolventaction}
\rho_{E,\chi}(ax)=\widehat a(\chi)\rho_{E,\chi}(x).
\end{equation}

\begin{lemma}[Character resolvent of a Gaussian period]\label{lem:fourierperiod}
Let \(m\in\cD(f)\) and \(\chi\in X\).  If \(f_\chi\mid m\), then
\begin{equation}\label{eq:echieta}
\rho_{E,\chi}(\eta_{E,m})
=\frac{1}{d_{E,m}}\tau_m(\chi^{-1});
\end{equation}
if \(f_\chi\nmid m\), the resolvent is zero.
Moreover, the resolvent is nonzero precisely when \(\chi\in\Phi_m\).
In particular,
\begin{equation}\label{eq:echitheta}
\rho_{E,\chi}(\Theta_E)
=\frac{1}{d_{E,m(\chi)}}\tau_{m(\chi)}(\chi^{-1}).
\end{equation}
\end{lemma}

\begin{proof}
Let \(N=\operatorname{Gal}(E/E_m)\) and
\(H_m=\operatorname{Gal}(\Q(\zeta_m)/E_m)\).
By Lemma~\ref{lem:dmcharacters}, \(\chi\) factors through
\(\operatorname{Gal}(E_m/\Q)\) exactly when \(f_\chi\mid m\).
If \(f_\chi\nmid m\), then \(\chi\) is nontrivial on \(N\).
In the sum \eqref{eq:resolvent}, group the terms into cosets of \(N\).
Since \(\eta_{E,m}\in E_m\) and
\(\sum_{h\in N}\chi(h)^{-1}=0\), the resolvent vanishes.

Assume \(f_\chi\mid m\).  Then \(\chi\) factors through
\(\operatorname{Gal}(E_m/\Q)\), each element of which has \(n/d_{E,m}\)
lifts to \(G\).  Expanding
\(\eta_{E,m}=\sum_{h\in H_m}\zeta_m^h\) therefore gives
\[
\rho_{E,\chi}(\eta_{E,m})
=\frac{1}{d_{E,m}}
\sum_{a\in(\Z/m\Z)^\times}\chi(a)^{-1}\zeta_m^a
=\frac{1}{d_{E,m}}\tau_m(\chi^{-1}).
\]
Lemmas~\ref{lem:imprimitivegauss} and \ref{lem:phimcriterion} show that
this is nonzero exactly for \(\chi\in\Phi_m\).
There is exactly one such index \(m(\chi)\), so summing over \(m\) proves
\eqref{eq:echitheta}.  
\end{proof}

\begin{proposition}[Fourier formula for the integral trace]\label{prop:traceformula}
Let \(a,b\in\cA_E\).  Then
\begin{equation}\label{eq:traceformula}
\operatorname{Tr}_{E/\Q}\bigl((a\Theta_E)(b\Theta_E)\bigr)
=
|G|\sum_{\chi\in X}
\widehat a(\chi)\widehat b(\chi^{-1})
\frac{\chi(-1)f_\chi}{d_{E,m(\chi)}^2}.
\end{equation}
\end{proposition}

\begin{proof}
Character orthogonality gives Fourier inversion in the form
\[
g(x)=\sum_{\chi\in X}\chi(g)\rho_{E,\chi}(x)
\qquad (g\in G,\ x\in E).
\]
Substituting \eqref{eq:resolvent}, the coefficient of \(h(x)\) is
\(n^{-1}\sum_{\chi\in X}\chi(gh^{-1})\).
This is \(1\) for \(h=g\) and \(0\) otherwise.  It follows that
\begin{equation}\label{eq:traceeigen}
\begin{aligned}
\operatorname{Tr}_{E/\Q}(xy)
&=\sum_{g\in G}g(x)g(y)\\
&=\sum_{\chi,\psi\in X}
\rho_{E,\chi}(x)\rho_{E,\psi}(y)
\sum_{g\in G}(\chi\psi)(g)\\
&=n\sum_{\chi\in X}
\rho_{E,\chi}(x)\rho_{E,\chi^{-1}}(y).
\end{aligned}
\end{equation}
All products on the right are products of complex numbers under the fixed
embedding, and the left side is the usual rational field trace.

Now take \(x=a\Theta_E\) and \(y=b\Theta_E\).  Since
\(f_{\chi^{-1}}=f_\chi\), we also have \(m(\chi^{-1})=m(\chi)\).
Using \eqref{eq:resolventaction} and \eqref{eq:echitheta}, the corresponding
summand in \eqref{eq:traceeigen} is
\[
\widehat a(\chi)\widehat b(\chi^{-1})
\frac{\tau_{m(\chi)}(\chi^{-1})\tau_{m(\chi)}(\chi)}
{d_{E,m(\chi)}^2}.
\]
Lemma~\ref{lem:gaussproduct} proves \eqref{eq:traceformula}.  
\end{proof}

\section{Comparison of integral trace lattices}\label{sec:comparison}

We now show directly that preservation of the conductor and parity of every
character gives an explicit isometry of integral trace lattices.

\begin{theorem}[Comparison from conductors and parities]\label{thm:comparison}
Let \(E\) and \(F\) be finite abelian extensions of \(\Q\).  Suppose that there is
a group isomorphism
\[
\iota\colon X_E\xrightarrow{\sim}X_F
\]
preserving the conductor and parity of every character.  Then \(T_E\cong T_F\).
The isometry can be chosen equivariant for the group isomorphism
\(J\colon G_E\xrightarrow{\sim}G_F\) characterized by
\begin{equation}\label{eq:dualcompatibility}
\iota(\chi)(J(g))=\chi(g)
\qquad (\chi\in X_E,\ g\in G_E).
\end{equation}
\end{theorem}

\begin{proof}
Finite abelian character duality gives the unique isomorphism \(J\) in
\eqref{eq:dualcompatibility}.  Thus \(J\) is the inverse of the map
\(G_F\to G_E\) dual to \(\iota\).
Extend \(J\) linearly to a \(\Q\)-algebra isomorphism
\[
J_*\colon\Q[G_E]\xrightarrow{\sim}\Q[G_F],
\qquad
J_*\left(\sum_g a_g g\right)=\sum_g a_gJ(g).
\]
Since \(J\) is a group isomorphism, this extension sends \(\Z[G_E]\)
onto \(\Z[G_F]\) and respects multiplication.

We use field subscripts or superscripts to distinguish the sets \(\Phi_m\)
and the corresponding idempotents for \(E\) and \(F\).
Conductor preservation gives a common conductor
\(f=f_E=f_F\), a common set \(\cD(f)\), and
\(\iota(\Phi_{E,m})=\Phi_{F,m}\).
After extending \(J_*\) to complex coefficients, \eqref{eq:dualcompatibility}
gives
\[
\begin{aligned}
J_*(e_\chi^E)
&=\frac1{|G_E|}\sum_{g\in G_E}\chi(g)^{-1}J(g)\\
&=\frac1{|G_F|}\sum_{h\in G_F}\iota(\chi)(h)^{-1}h
=e_{\iota(\chi)}^F.
\end{aligned}
\]
Consequently \(J_*(e_{E,m})=e_{F,m}\).
Theorem~\ref{thm:leopoldt} now implies
\[
J_*(\cA_E)=\cA_F.
\]
The same theorem says that each element of \(\cO_E\) has a unique
representation \(a\Theta_E\) with \(a\in\cA_E\), and similarly for \(F\).
Therefore
\begin{equation}\label{eq:comparisonmap}
U\colon\cO_E\longrightarrow\cO_F,
\qquad U(a\Theta_E)=J_*(a)\Theta_F,
\end{equation}
is a well-defined \(\Z\)-linear isomorphism.  Its inverse is
\(b\Theta_F\mapsto J_*^{-1}(b)\Theta_E\).
For \(g\in G_E\) and \(a\in\cA_E\), it also satisfies
\[
U\bigl(g(a\Theta_E)\bigr)
=J_*(ga)\Theta_F
=J(g)J_*(a)\Theta_F
=J(g)U(a\Theta_E).
\]

It remains to verify preservation of the pairing.
Lemma~\ref{lem:dmcharacters} and conductor preservation give
\[
d_{E,m}
=\#\{\chi\in X_E:f_\chi\mid m\}
=\#\{\psi\in X_F:f_\psi\mid m\}
=d_{F,m}.
\]
Moreover, \(m_F(\iota(\chi))=m_E(\chi)\), and the Fourier transforms satisfy
\[
\widehat{J_*(a)}(\iota(\chi))
=\sum_{g\in G_E}a_g\iota(\chi)(J(g))
=\widehat a(\chi).
\]
Finally, \(\iota(\chi^{-1})=\iota(\chi)^{-1}\) and
\(\iota(\chi)(-1)f_{\iota(\chi)}=\chi(-1)f_\chi\).
Thus, relabelling the sum in \eqref{eq:traceformula} by
\(\psi=\iota(\chi)\) gives, for every \(a,b\in\cA_E\),
\[
\operatorname{Tr}_{F/\Q}
\bigl((J_*(a)\Theta_F)(J_*(b)\Theta_F)\bigr)
=
\operatorname{Tr}_{E/\Q}\bigl((a\Theta_E)(b\Theta_E)\bigr).
\]
By \eqref{eq:comparisonmap}, this is precisely
\(B_F(Ux,Uy)=B_E(x,y)\) for \(x,y\in\cO_E\).
Hence \(U\) is the required equivariant isometry.  
\end{proof}

\begin{remark}\label{rem:strongerabelian}
Theorem~\ref{thm:comparison} does not require the Galois group to be cyclic.  The cyclic hypothesis enters only in Theorem~\ref{thm:characterdata}, where the degree, discriminant, and signature are shown to recover a conductor- and parity-preserving identification of the character groups.
\end{remark}

\section{The main theorem and its consequences}\label{sec:mainproof}

We now deduce the cyclic classification and record its wild and tame
consequences together with the role of parity and the associated order.

\begin{proof}[Proof of Theorem~\ref{thm:main}]
Suppose first that \(T_K\cong T_L\).  Proposition~\ref{prop:necessary} gives
\[
[K:\Q]=[L:\Q],
\qquad
\operatorname{sig}(K)=\operatorname{sig}(L),
\qquad
\operatorname{Disc}(K)=\operatorname{Disc}(L).
\]

Conversely, suppose these three arithmetic invariants agree.  Theorem~\ref{thm:characterdata} gives an isomorphism
\[
\iota\colon X_K\xrightarrow{\sim}X_L
\]
preserving the conductor and parity of every character.  Theorem~\ref{thm:comparison} then gives \(T_K\cong T_L\).  
\end{proof}

\begin{corollary}[Wildly ramified cyclic fields]\label{cor:wild}
Let \(K\) and \(L\) be wildly ramified cyclic number fields of the same degree, signature, and discriminant.  Then \(T_K\cong T_L\).
\end{corollary}

\begin{proof}
This follows directly from Theorem~\ref{thm:main}; no tameness assumption
occurs in its proof.  
\end{proof}

\begin{corollary}[Tame classification]\label{cor:tame}
Let \(K\) and \(L\) be tame cyclic number fields of the same degree and discriminant.  Then \(T_K\cong T_L\).
\end{corollary}

\begin{proof}
A tame abelian extension of \(\Q\) is unramified at \(2\): the order
\(|Y_{K,2}|\) is both a power of \(2\) and prime to \(2\), so \(Y_{K,2}=1\).
Every ramified odd-prime local subgroup is recovered from the discriminant by
Proposition~\ref{prop:odddetermination}.  Hence the parity of a global generator,
and thus the signature, is determined by the degree and discriminant.
Theorem~\ref{thm:main} applies.  This recovers the principal classification
theorem of \cite{BolanosMantilla}.  
\end{proof}

\begin{remark}[Where wild ramification enters]\label{rem:wild}
In the tame case, all character conductors are squarefree and Leopoldt's associated order is \(\Z[G]\), so the argument can be implemented using a normal integral basis.  In the wild case, \(\cO_K\) need not be free over \(\Z[G]\).  The role of Leopoldt's theorem is precisely to replace \(\Z[G]\) by the associated order \(\cA_K\), over which \(\cO_K\) is free of rank one.  Formula \eqref{eq:traceformula} shows that the additional data appearing in the wild case are already encoded by the conductors of the characters.
\end{remark}

\begin{example}[Equal signed discriminants and different signatures]
\label{ex:quartic}
Let
\[
K_+=\Q(\alpha),\quad \alpha=\sqrt{2+\sqrt2},
\qquad
K_-=\Q(\beta),\quad \beta=\sqrt{-2-\sqrt2}.
\]
Both fields are cyclic quartic fields of conductor \(16\) and signed
discriminant \(2^{11}\).  Their integral trace lattices are not isometric.
\end{example}

\begin{proof}
Put \(\zeta=\exp(2\pi i/16)\), and choose
\(\alpha=\zeta+\zeta^{-1}\) and \(\beta=\zeta^3-\zeta^{-3}\).
Their respective minimal polynomials are
\[
P_+(t)=t^4-4t^2+2,
\qquad P_-(t)=t^4+4t^2+2,
\]
which are Eisenstein at \(2\).
In \(\operatorname{Gal}(\Q(\zeta)/\Q)=(\Z/16\Z)^\times\), the elements
fixing \(K_+\) form \(\{1,15\}\), whereas those fixing \(K_-\) form \(\{1,7\}\).
Indeed, the indicated nonidentity automorphism fixes the displayed
generator in each case; its stabilizer has order \(2\) because the generator
has degree \(4\).
Since \((\Z/16\Z)^\times=\langle-1\rangle\times\langle5\rangle\),
the corresponding character groups are
\[
X_{K_+}=\langle\psi_2\rangle,
\qquad X_{K_-}=\langle\varepsilon\psi_2\rangle.
\]
They are cyclic of order \(4\).  In each group the character conductors are
\(1,8,16,16\), so the field conductor is \(16\) and the absolute discriminant
is \(8\cdot16^2=2^{11}\) by \eqref{eq:conductor-discriminant}.
The roots of \(P_+\) are all real, while the roots of \(P_-\) are all purely
imaginary.  Thus the field signatures are \((4,0)\) and \((0,2)\), respectively,
and both signed discriminants are positive.
Proposition~\ref{prop:necessary} gives real trace-form signatures \((4,0)\)
and \((2,2)\), so the trace lattices cannot be isometric.

For completeness, \(\operatorname{disc}(P_+)=\operatorname{disc}(P_-)=2^{11}\).
Comparison with the field discriminants shows that the respective power
bases are integral bases.  In those bases the trace Gram matrices are
\[
\begin{pmatrix}
4&0&8&0\\
0&8&0&24\\
8&0&24&0\\
0&24&0&80
\end{pmatrix}
\quad\text{and}\quad
\begin{pmatrix}
4&0&-8&0\\
0&-8&0&24\\
-8&0&24&0\\
0&24&0&-80
\end{pmatrix},
\]
each of determinant \(2^{11}\).  These matrices follow from Newton's
identities applied to \(P_+\) and \(P_-\).  
\end{proof}

\begin{remark}[Necessity of the parity hypothesis]\label{rem:paritynecessary}
The isomorphism sending \(\psi_2^j\) to \((\varepsilon\psi_2)^j\) preserves the
conductor of every character in Example~\ref{ex:quartic}, but does not
preserve the parity of the order-\(4\) characters.  Thus the parity assumption
cannot be omitted from Theorem~\ref{thm:comparison}.  This is a failure of
the universal comparison statement without that assumption, not a claim
that every isometry of abelian trace lattices must arise from a
character-group isomorphism.  The same example shows that the field-signature
hypothesis in Theorem~\ref{thm:main} cannot be omitted, even when one fixes
the signed discriminant.
\end{remark}

\section*{Declaration of AI use}
ChatGPT (OpenAI) was used for assistance with literature searches,
exploration and review of mathematical arguments, and manuscript
preparation.  The author is responsible for the mathematical content,
accuracy of the references, and final version submitted for publication.

\end{document}